\documentclass{amsart}

\newcommand{\Real}{\mathrm{Re}}

\begin{document}

\title{An uncertainty principle related to the Euclidean motion group}
\author{Jens Gerlach Christensen}
\author{Henrik Schlichtkrull}
\address{Department of mathematics, University of Copenhagen}
\email{vepjan@math.ku.dk, schlicht@math.ku.dk}
\urladdr{http://www.math.ku.dk/\textasciitilde vepjan,
http://www.math.ku.dk/\textasciitilde schlicht}

\begin{abstract}
  We show that a well known uncertainty principle for 
  functions on the circle can be derived from an uncertainty 
  principle for the Euclidean motion group.
\end{abstract}

\maketitle

\section{Uncertainty principles related to Lie group representations}
Let $G$ be a Lie group with Lie algebra $\mathfrak{g}$, and let
$(\pi,H)$ be a unitary representation of $G$. Then each element
$X\in\mathfrak{g}$ generates a closed, skew-adjoint operator
$\pi(X)$ on $H$ by
\begin{equation*}
  \pi(X)f(s) = \lim_{t\to 0} \frac{\pi(\exp tX)v - v}{t},
\end{equation*}
with domain $D(\pi(X))$ consisting of
all $v\in H$ for which the limit exists.

The uncertainty principle related to $\pi$ 
tells that for operators generated by
$X$,$Y$ and $[X,Y]$ the following holds
\begin{equation*}
    \|\pi(X)f\|\|\pi(Y)x\| 
      \geq \frac{1}{2}|\langle \pi([X,Y])x,x\rangle |
\end{equation*}
for all $x\in D(\pi(X))\cap D(\pi(Y))\cap D(\pi([X,Y]))$.

We would like to advocate this as a natural way to achieve
uncertainty principles. It was first proposed in Kraus \cite{kraus}
for Lie groups of dimension $\leq 3$, and the dimension constraint was 
recently removed by the first author in \cite{uplie}.
For example, the classical Heisenberg uncertainty principle for
functions on $\mathbb{R}^n$
is easily derived in this way from the Schr\"odinger representation
of the Heisenberg group, see \cite{folsit}, p.\ 212. 

It is the purpose of this note to point out that the uncertainty
principle for the circle, which was motivated by Breitenberger
\cite{breit} and further discussed in \cite{narcward}, \cite{prequak},
\cite{pqrs}, \cite{rauhut},  \cite{margit}, \cite{selig}
is obtained similarly
from the principal series representation of the Euclidean motion 
group of $\mathbb{R}^2$.

\section{The Euclidean motion group and a unitary representation}
Let $G$ be the Euclidean motion group
\begin{equation*}
  G = \left\{ 
      (r,z) = 
      \begin{pmatrix}
        e^{ir} & z \\
        0      & 1 
      \end{pmatrix}
      \Big| r\in\mathbb{R}, z\in\mathbb{C}
      \right\}
\end{equation*}
Its Lie algebra is
\begin{equation*}
  \mathfrak{g} = 
      \left\{ 
      \begin{pmatrix}
        ir & z \\
        0  & 0 
      \end{pmatrix}
      \Big| r\in\mathbb{R}, z\in\mathbb{C}
      \right\}
\end{equation*}
Let $H$ be the Hilbert space $H = L_2(\mathbb{T}),$ of
square integrable functions on the circle
$\mathbb{T}=\{s\in\mathbb{C}|\mid |s|=1\}$.
As in \cite[Chapter IV]{sugiura}
the following defines a unitary representation of $G$ on $H$:
\begin{equation*}
  \pi_a(r,z)f(s) = e^{i(z,sa)}f(e^{-ir}s), \quad (s\in\mathbb{T})
\end{equation*}
where $a\in\mathbb{C}$ and $(x,y)=\Real (x\overline{y}).$ 
For simplicity we assume in the following that $a=1$, which
is sufficient for our purpose. The representation $\pi_1$ will be
denoted $\pi$.
%
\section{Operators generated from the representation}

We now generate three operators from elements of the Lie algebra $\mathfrak{g}$.
Let $X,Y_1,Y_2\in\mathfrak{g}$ be
\begin{equation*}
  X =  
      \begin{pmatrix}
        i & 0 \\
        0 & 0 
      \end{pmatrix}
  \qquad\text{,}\qquad 
  Y_1 = 
      \begin{pmatrix}
        0 & 1 \\
        0 & 0 
      \end{pmatrix}
  \qquad\text{and}\qquad 
  Y_2 = 
      \begin{pmatrix}
        0 & i \\
        0 & 0 
      \end{pmatrix}
\end{equation*}
then
\begin{equation*}
  \exp(tX) =
      \begin{pmatrix}
        e^{it} & 0 \\
        0      & 1 
      \end{pmatrix}
  \quad\text{,}\quad 
  \exp(tY_1) = 
      \begin{pmatrix}
        1 & t \\
        0 & 1 
      \end{pmatrix}
  \quad\text{and}\quad 
  \exp(tY_2) = 
      \begin{pmatrix}
        1 & it \\
        0 & 1 
      \end{pmatrix}
\end{equation*}
and we then get
\begin{equation*}
  \pi(X)f(s)
  = \lim_{t\to 0} \frac{\pi(t,0)f(s) - f(s)}{t}
  = \lim_{t\to 0} \frac{f(e^{-it}s) - f(s)}{t}
  = -f'(s),
\end{equation*}
where $f'(s)=\frac{d}{dt}f(e^{it}s).$
This operator has domain 
\begin{equation}\label{domain}
  \{ f\in L_2(\mathbb{T}) | 
          \text{$t\mapsto f(e^{it})$ absolutely continuous with }
f'\in L_2(\mathbb{T}) \} 
\end{equation}
Also 
\begin{equation*}
  \pi(Y_1)f(s)
  = \lim_{t\to 0} \frac{\pi(0,t)f(s) - f(s)}{t}
  = \lim_{t\to 0} \frac{e^{i(t,s)}f(s) - f(s)}{t}
  = i\cos(\theta)f(s)
\end{equation*}
and
\begin{equation*}
  \pi(Y_2)f(s)
  = \lim_{t\to 0} \frac{\pi(0,it)f(s) - f(s)}{t}
  = \lim_{t\to 0} \frac{e^{i(it,s)}f(s) - f(s)}{t}
  = i\sin(\theta)f(s)
\end{equation*}
when $s=e^{i\theta}.$ Both of these operators are defined on 
$H$. 

\section{The uncertainty principle}
Since 
\begin{equation*}
 [X,Y_1]=Y_2,\quad [X,Y_2]=-Y_1
\end{equation*}
the uncertainty principle gives
\begin{equation} \label{rauhut1}
  |\langle \pi(Y_1)x,x\rangle| =
  |\langle \pi([X,Y_2])x,x\rangle| \leq 2 
\| \pi(X)x\|\|\pi(Y_2)x\|
\end{equation}
and
\begin{equation} \label{rauhut2}
  |\langle \pi(Y_2)x,x\rangle| =
  |\langle \pi([X,Y_1])x,x\rangle| \leq 2 
\| \pi(X)x\|\|\pi(Y_1)x\|.
\end{equation}

Let $T=i\pi(X)$ denote the operator $Tf=-if'$ with domain
(\ref{domain}),
and let $S_1=i\pi(Y_1)$, $S_2=i\pi(Y_2)$, then
$T$, $S_1$ and $S_2$ are selfadjoint, and the unitary
operator $S=S_1+iS_2$ 
is given by $Sf(e^{i\theta})=-e^{i\theta}f(e^{i\theta})$
with $D(S)=H$.

Now
\begin{equation*}
\| S x\|^2 = \| S_1x\|^2 + \| S_2 x\|^2 = \|x\|^2
\end{equation*}
and 
\begin{equation*}
  |\langle Sx,x\rangle |^2 = 
  |\langle S_1x,x\rangle |^2 + |\langle S_2x,x\rangle |^2  
\end{equation*}
so the uncertainty principles (\ref{rauhut1}) and (\ref{rauhut2}) 
give
\begin{equation*}
  |\langle Sx,x\rangle|^2 
  \leq 4\| x\|^2 \|Tx\|^2.
\end{equation*}

As explained in \cite[Theorem 2.3.3]{rauhut} 
this is exactly the uncertainty principle of Breitenberger.
Notice that in this case we have a $3$-dimensional Lie group so that the
example is covered by \cite{kraus}.

\end{document}